\documentclass[12p,preprint]{elsarticle}

\usepackage{amsmath}
\usepackage{amsfonts}
\usepackage{amssymb,latexsym}
\usepackage{amsthm}
\usepackage{graphicx}

\usepackage{arcs}
\usepackage{booktabs}
\usepackage{caption}
\usepackage{color}
\usepackage{hyperref}
\usepackage{multirow}
\usepackage{newlfont}
\usepackage{rotating}
\usepackage{subfigure}
\usepackage{verbatim}
\usepackage{xcolor} 
\usepackage{units}
\usepackage{bm}
\usepackage[english]{babel}
\usepackage[utf8]{inputenc}
\usepackage{algorithm}
\usepackage[noend]{algpseudocode}
\usepackage{cancel}
\usepackage[normalem]{ulem}

\newtheorem{Proposition}{Proposition}
\newtheorem{Lemma}{Lemma}
\newtheorem{Remark}{Remark}

\newcommand{\NN}{\mathbb{N}}
\newcommand{\ZZ}{\mathbb{Z}}
\newcommand{\QQ}{\mathbb{Q}}
\newcommand{\RR}{\mathbb{R}}
\newcommand{\CC}{\mathbb{C}}
\newcommand{\II}{\mathbb{I}}
\newcommand{\PP}{\mathbb{P}}

\DeclareMathOperator{\csch}{csch}

\newcommand{\red}[1]{{\color{red} #1}}
\newcommand{\blue}[1]{{\color{blue} #1}}

\begin{document}

\begin{frontmatter}

\title{Chebyshev admissible meshes and Lebesgue constants\\ of complex polynomial projections}



\address[address-PL]{Jagiellonian University, Poland}
\address[address-PD]{University of Padova, Italy}

\author[address-PL]{Leokadia Białas-Cież}
\ead{leokadia.bialas-ciez@uj.edu.pl}

\author[address-PL]{Dimitri Jordan Kenne}
\ead{dimitri.kenne@doctoral.uj.edu.pl}

\author[address-PD]{Alvise Sommariva}
\ead{alvise@math.unipd.it}

\author[address-PD]{Marco Vianello}
\ead{marcov@math.unipd.it}


\begin{abstract}
We construct admissible polynomial meshes on piecewise polynomial or trigonometric curves of the complex plane, by mapping univariate Chebyshev points. Such meshes can be used for polynomial least-squares, for the extraction of Fekete-like and Leja-like interpolation sets, and also for the evaluation of their Lebesgue constants.
\end{abstract}



\end{frontmatter}

\section{Complex Chebyshev-like polynomial meshes}

Starting from the seminal paper by Calvi and Levenberg \cite{CL08}, the notion of {\em polynomial (admissible) mesh} has been emerging in the last years as a fundamental theoretical and computational tool in polynomial approximation. In the present paper we focus on the univariate complex case. We recall that an admissible polynomial mesh of a polynomially determining compact set $K\subset \mathbb{C}$ (i.e., a polynomial vanishing on $K$ vanishes everywhere on $\mathbb{C}$), is a sequence of finite norming subsets $Z_n\subset K$ such that 
\begin{equation} \label{am}
\|p\|_K\leq c \|p\|_{Z_n}\;,\;\;\forall p\in \mathbb{P}_n(\mathbb{C})\;,
\end{equation}
where $\|\cdot\|_Y$ denotes the sup-norm on a continuous or discrete compact set $Y$, $card(Z_n) = O(n^\alpha)$, $\alpha\geq 1$, and $c$ is a constant independent
of $n$. The fact that $card(Z_n)\geq dim(\mathbb{P}_n(\mathbb{C}))=n+1$ necessarily holds, since each $Z_n$ is $\mathbb{P}_n(\mathbb{C})$-determining. Such a mesh is termed {\em optimal} when $\alpha = 1$.

To give only a flavour of the topic, we recall that polynomial meshes are invariant by affine transformations, are stable under small perturbations, and can be 
assembled by finite union, finite product and algebraic transformations, 
starting from known instances. Moreover, admissible meshes can be conveniently used for least-square approximation, and contain extremal sets for interpolation of Fekete and Leja type, that can be computed by greedy algorithms; cf., e.g., \cite{BCC12,BCLSV11,CL08, SV09}. 
 
Existence of admissible meshes with $\mathcal{O}(n^2)$ cardinality has been proved on any connected compact set of $\mathbb{C}$ whose boundary is a $C^1$ parametric curve with bounded tangent vectors, while optimal admissible meshes are known in special instances; cf. \cite{BCC12}. The following Proposition and Remark show how to construct optimal admissible meshes 
of Chebyshev type on a wide class of complex curves and domains. To this purpose, we need a basic Lemma.

\begin{Lemma}
Let $\phi(t)$, $t\in [a,b]$, be an algebraic or trigonometric polynomial with complex coefficients, of degree not exceeding $\nu$ (with $b-a\leq 2\pi$ in the trigonometric case). Denote by $\mathcal{T}_N$ the set of $N$ Chebyshev zeros in $(-1,1)$, 
$\cos((2j-1)\pi/(2N))$, $1\leq j\leq N$, or the set of $N+1$ Chebyshev extrema 
in $[-1,1]$, $\cos(j\pi/N)$, $0\leq j\leq N$.

Consider the points 
\begin{equation} \label{mesh}
X_\nu^m=\sigma(\mathcal{T}_N) \subset [a,b]
\end{equation}
where
\begin{equation} \label{alg}
N=m\nu\;,\;\;\sigma(u)=\frac{b-a}{2}\,u+\frac{b+a}{2}\;,
\;\;u\in [-1,1]\;,
\end{equation}
in the algebraic case, and 
\begin{equation} \label{trig}
N=2m\nu\;,\;\;\sigma(u)=2\arcsin\left(u\sin\left(\frac{b-a}{4}\right)\right)+\frac{b+a}{2}\;,\;\;u\in [-1,1]\;,
\end{equation}
in the trigonometric case.

Then the following inequality holds for every $\nu\geq 1$, $m>1$
\begin{equation} \label{pmesh}
\|\phi\|_{[a,b]}\leq c_m\|\phi\|_{X_\nu^m}\;,\;\;c_m:=\frac{1}{\cos(\pi/(2m))}\;.
\end{equation}
\end{Lemma}
\vskip0.3cm
\noindent
{\em Proof.} 
In the {\em real} algebraic case, (\ref{pmesh}) is a well-known polynomial inequality originally proved by Ehlich and Zeller \cite{EZ64}. Interestingly, a proof can be given also by the notion of {\em Dubiner distance} in $[a,b]$, which is tailored to polynomial spaces; cf., e.g., \cite{B17,PV18}. Indeed, in \cite{V18-2} such a notion has been extended in the subperiodic trigonometric case, i.e. to real trigonometric polynomials on subintervals of the period, namely on $[a,b]$ with $b-a\leq 2\pi$. In such a way, inequality (\ref{pmesh})  has been proved for {\em real} trigonometric polynomials.  

We now show how to extend such inequality to  algebraic and trigonometric polynomials of a real variable with {\em complex} coefficients. 
Take $t^\ast\in [a,b]$ such that $|\phi(t^\ast)|=\|\phi\|_{[a,b]}$. We can assume that $\phi(t^\ast)\neq 0$, since (\ref{pmesh}) trivially holds for $\phi\equiv 0$. Define the complex number 
$u=\overline{\phi(t^\ast)}/|\phi(t^\ast)|$ (which lies on the unit circle), and observe that $u\phi(t^\ast)=|\phi(t^\ast)|^2/|\phi(t^\ast)|=|\phi(t^\ast)|$. 

Now, consider $\psi(t)=u\phi(t)$; clearly, $|\psi(t)|=|\phi(t)|$, and $Im(\psi(t^\ast))=0$, since $\psi(t^\ast)=|\phi(t^\ast)|$ is real. Since $Re(\psi(t))$ is a {\em real} algebraic or trigonometric polynomial, we can write the chain of inequalities
\begin{eqnarray}
\|\phi\|_{[a,b]}&=&|\phi(t^\ast)|=Re(\psi(t^\ast)) \leq \|Re(\psi)\|_{[a,b]}
\leq c_m\|Re(\psi)\|_{X_\nu^m} \nonumber \\
&\leq& c_m\|\psi\|_{X_\nu^m}=c_m\|\phi\|_{X_\nu^m}\;.\hspace{4cm} \square \nonumber
\end{eqnarray}

\begin{Proposition}
Let $\Gamma$ be (the image of) a complex parametric curve $z(t)$, $t\in [a,b]$, where $z(t)$ is an  algebraic or trigonometric polynomial of degree $k\geq 1$ (with $b-a\leq 2\pi$ in the trigonometric case). 

Then the sequence $Z_n^m(k)=z(X_{nk}^m)$, cf. (\ref{mesh})-(\ref{trig}), forms an (optimal) admissible polynomial mesh for $\Gamma$, since the following polynomial inequality holds for every 
$p\in \mathbb{P}_n(\mathbb{C})$, $n\geq 1$, $m>1$
\begin{equation} \label{pmesh1}
\|p\|_\Gamma\leq c_m\|p\|_{Z_n^m(k)}\;.
\end{equation}
\end{Proposition}
\vskip0.3cm
\noindent
{\em Proof.} Consider the function composition $\phi(t)=p(z(t))$,  which clearly is an algebraic or trigonometric polynomial on $[a,b]$ with {\em complex coefficients}, of degree at most $\nu=nk$. The result is an immediate consequence of Lemma 1, by observing that
$$
\|p\|_\Gamma=\|\phi\|_{[a,b]}\leq c_m\|\phi\|_{X_{nk}^m}=c_m\|p\|_{Z_n^m(k)}\;. \hspace{0.3cm} \square
$$

\begin{Remark}
Let $\Gamma=\bigcup_{j=1}^s{\Gamma_j}$ be union of parametric algebraic or trigonometric arcs $\Gamma_j$ of degree $k_j$ on $[a_j,b_j]$, $1\leq j\leq s$. Then for every  
$p\in \mathbb{P}_n(\mathbb{C})$, $n\geq 1$, $m>1$ 
\begin{equation} \label{union}
\|p\|_\Gamma\leq c_m\|p\|_{Z_n^m}\;,\;\;Z_n^m=\bigcup_{j=1}^s{Z_n^m(k_j)}\;,
\end{equation}
i.e. $Z_n^m$ is an optimal admissible mesh for $\Gamma$, by the finite union property of admissible meshes, cf. e.g. \cite[Lemma 4]{CL08}. On the other hand, such $Z_n^m$ is an admissible mesh also for any compact set $K\subset \mathbb{C}$ having outer boundary lying on $\Gamma$, with $\Gamma$ contained in $K$, say $\partial K_{out}\subseteq \Gamma\subseteq K$, by the maximum modulus principle applied to polynomials (we recall that the outer boundary is the boundary of the unbounded connected component of $\mathbb{C}\setminus K$). Notice that such a class is very wide: it includes linear polygons, as well as curvilinear polygons with boundary tracked by splines, or by arcs like $r(\theta)(\cos(\theta),\sin(\theta))$ in polar coordinates with $r(\theta)$ a trigonometric polynomial. See the Figures below for some illustrative examples.
\end{Remark}

The following Proposition shows that suitable admissible meshes of the form (\ref{pmesh1}) can be conveniently used to evaluate Lebesgue constants, with rigorous error bounds. Again, we begin with a basic Lemma. The result is well-known for interpolation operators (cf. e.g. \cite{T23}), nevertheless we prefer to prove it here for more general projection operators (which include for example also least-square approximations). Below by $C(K)$ we denote as usual the space of continuous functions on the compact set $K\subset \mathbb{C}$.

\begin{Lemma}
Let $K\subset \mathbb{C}$ be a compact set and  $L_n:C(K)\to \mathbb{P}_n(\mathbb{C})$ a linear projection operator such that 
\begin{equation} \label{Ln}
L_nf(z)=\sum_{j=1}^M{f(\xi_j)\,\phi_j(z)}\;,
\end{equation}
where $\Xi=\{\xi_j\}\subset K$ and $\{\phi_j\}$ is a set of generators of $\mathbb{P}_n(\mathbb{C})$.  Moreover, let 
$$
\lambda_n(z)=\sum_{j=1}^M{|\phi_j(z)|}
$$
be the ``Lebesgue function'' of $L_n$.  

Then the ``Lebesgue constant'' of $L_n$, that is its uniform norm, is equal to the sup-norm of the Lebesgue function on $K$
 $$
 \|L_n\|=\sup_{f\neq 0}\frac{\|L_nf\|_K}{\|f\|_K}=\|\lambda_n\|_K=\|\lambda_n\|_{ {\partial}K_{out}}\;.
 $$
\end{Lemma}
\vskip0.3cm
\noindent
{\em Proof.} Inequality $\|L_n\|\leq \|\lambda_n\|_K$ is immediate, since $$|L_nf(z)|\leq \sum_{j=1}^M{|f(\xi_j)|\,|\phi_j(z)}|\leq 
\|f\|_\Xi\lambda_n(z)\leq \|f\|_K\lambda_n(z)\;.$$
Let $z^\ast\in K$ such that $\|\lambda_n\|_K=|\lambda_n(z^\ast)|$. Now, the point is to find a continuous function $f^\ast$ on $K$ such that $f^\ast(\xi_j)=u_j=|\phi_j(z^\ast)|/\phi_j(z^\ast)$ for all $j$ such that $\phi_j(z^\ast)\neq 0$, and $\|f^\ast\|_K=1$. To this purpose, since $|u_j|=1$
let us write $u_j=e^{i\theta_j}$, where $\theta_j\in [0,2\pi)$, and define a function $g:\{\xi_j\}\to [0,2\pi)$ such that $g(\xi_j)=\theta_j$. By a deep topological result, the celebrated Tietze extension theorem (cf. e.g. \cite[Ch.7, Thm.5.1]{D66}), 
since $g$ is trivially continuous on the closed discrete subset $\{\xi_j\}$, there exists an extension $\tilde{g}\in C(K)$ taking values in $[0,2\pi)$. 
Then, $f^\ast(z)=e^{i\tilde{g}(z)}\in C(K)$ is the required function, because $f^\ast(\xi_j)=e^{i\tilde{g}(\xi_j)} =e^{i\theta_j}=u_j$ and $|f^\ast(z)|\equiv 1$. To prove that $\|\lambda_n\|_K=\|\lambda_n\|_{\partial K_{out}}$, we can clearly restrict to compact domains (the closure of bounded connected open sets). Then we can apply the maximum principle for subharmonic functions to $\lambda_n$, since each $|\phi_j|$ is subharmonic being the modulus of an (entire) holomorphic function and the sum of subharmonic functions is subharmonic; cf. e.g. \cite[\S7.7]{K08}.
\hspace{0.3cm} $\square$
\vskip0.3cm

\begin{Remark}
Existence of a continuous function $f^\ast$ as in the proof above can also be proved by Dugundji's version of Tietze extension theorem, which in its general formulation concerns extension of continuous functions defined on closed subsets of metric spaces and taking values in locally convex topological vector spaces; cf. \cite[Ch.9,Thm.6.1]{D66}. Applied to the present context, it simply says that defining $f$ such that $f(\xi_j)=u_j$, there exists an extension $f^\ast\in C(K)$ such that $f^\ast(K)\subset convhull(\{f(\xi_j)\})$. Thus $\|f^\ast\|_K=1$, because the 
$u_j$ lie on the unit circle in $\mathbb{C}$ and hence their convex hull is a polygon lying in the unit disk.
\end{Remark}

\begin{Remark}
The structure of projection operators like (\ref{Ln}) includes interpolation operators at $n+1$ distinct nodes $\xi_1,\dots,\xi_{n+1}$, where, denoting by 
$V_n=[p_j(\xi_i)]$, $1\leq i,j \leq n+1$, the Vandermonde-like matrix in any fixed polynomial basis $span(p_1,\dots,p_{n+1})=\mathbb{P}_n(\mathbb{C})$, we have that 
\begin{equation} \label{interp}
\phi_j(z)=\ell_j(z)=\frac{det(V_n(\xi_1,\dots,\xi_{j-1},z,\xi_{j+1},\dots,\xi_n))}{det(V_n(\xi_1,\dots,\xi_{j-1},\xi_j,\xi_{j+1},\dots,\xi_n))}
=\prod_{k=1,k\neq j}^{n+1}{(z-\xi_k)/(\xi_j-\xi_k)}
\end{equation}
are the corresponding fundamental Lagrange polynomials. But also discrete weighted least-squares operators at $M>n+1$ nodes $\Xi=\{\xi_1,\dots,\xi_M\}$ with positive weights 
$W=\{w_1,\dots,w_M\}$ are included. Indeed, denoting by $\{\pi_k\}$, $1\leq k\leq n+1$, the orthonormal polynomials with respect to the corresponding discrete scalar product $(f,g)_{\ell^2 _W(\Xi)}=\sum_{j=1}^M{w_jf(\xi_j)\overline{g(\xi_j)}}$, we have that
\begin{equation} \label{ls}
L_nf(z)=\sum_{k=1}^{n+1}{(f,\pi_k)_{\ell^2 _W(\Xi)}\,\pi_k(z)}=
\sum_{j=1}^M{f(\xi_j)w_jK_n(z,\xi_j)}\;,
\end{equation}
i.e. $\phi_j(z)=w_jK_n(z,\xi_j)$, where $K_n(z,v)=\sum_{k=1}^{n+1}{\pi_k(z)\overline{\pi_k(v)}}$ is the reproducing kernel of the discrete scalar product. Notice that in this case (unless $M=n+1$ where least-squares approximation coincides with interpolation) the $\phi_j$ are linearly dependent, thus forming a set of generators of $\mathbb{P}_n(\mathbb{C})$.
\end{Remark}

\begin{Proposition}
Let the assumptions of Remark 1 be satisfied, and assume that $L_n:C(K)\to \mathbb{P}_n(\mathbb{C})$ is a linear projection operator as in Lemma 2. 

Then the following estimates hold for every $n\geq 1$, $m>1$ 
\begin{equation} \label{lebest}
\|\lambda_n\|_{Z_n^m}\leq \|L_n\|\leq 
c_m\|\lambda_n\|_{Z_n^m}\;,
\end{equation} 
and 
\begin{equation} \label{errest}
0\leq \|L_n\|
-\|\lambda_n\|_{Z_n^m}
\leq (c_m-1)\|L_n\|\;.
\end{equation}

\end{Proposition}
\vskip0.3cm
\noindent
{\em Proof.} Applying inequality (\ref{union}) to the polynomial $L_nf$, in view of the maximum modulus principle we get 
$$
\|L_nf\|_K=\|L_nf\|_{\partial K_{out}}\leq c_m\|L_nf\|_{Z_n^m}\;.
$$
On the other hand $|L_nf(z)|\leq \|f\|_\Xi\lambda_n(z)\leq \|f\|_{K} \lambda_n(z)$ and thus $\|L_nf\|_{K}\leq c_m\|f\|_{K} \|\lambda_n\|_{Z_n^m}$, from which we get immediately 
$$
\|L_n\|=\|\lambda_n\|_{K} 
\leq c_m\|\lambda_n\|_{Z_n^m}\;,
$$
and thus (\ref{lebest}) and (\ref{errest}), since $\|\lambda_n\|_{K}\geq \|\lambda_n\|_{Z_n^m}$  by inclusion.\hspace{0.3cm} $\square$
\vskip0.3cm

\begin{Remark}
Notice that $c_m\to 1$ and thus, if the sampling set $\Xi$ is independent of $m$, $\|\lambda_n\|_{Z_n^m}\to \|L_n\|$ as $m\to \infty$. On the other hand, in any case (\ref{errest}) gives the relative error estimate
\begin{equation} \label{asymptotic}
\frac{\|L_n\| -\|\lambda_n\|_{Z_n^m}}{\|L_n\|} 
\leq c_m-1=\frac{1-\cos(\pi/(2m))}{\cos(\pi/(2m))}\sim \frac{\pi^2}{8m^2}\approx \frac{1.23}{m^2}\;,
\end{equation}
that is a $\mathcal{O}(1/m^2)$ relative approximation of the Lebesgue constant by $\|\lambda_n\|_{Z_n^m}$. 
For example, with $m=4$ we already get the Lebesgue constant with a relative error less than $10\%$, thus correctly estimating its actual order of magnitude that is the relevant parameter in polynomial approximation. On the other hand (\ref{lebest}) gives also
the rigorous and {\em computable\/} absolute error estimate 
\begin{equation} \label{absest}
0\leq \|L_n\| 
-\|\lambda_n\|_{Z_n^m}\leq (c_m-1) \|\lambda_n\|_{Z_n^m}\;.
\end{equation}
\end{Remark}

\section{Numerical tests}

In this section we present several numerical tests (implemented in Matlab), concerning the use of complex Chebyshev-like polynomial meshes for interpolation and least-squares approximation, as well as for the evaluation of the corresponding Lebesgue constants. In all these applications we can work on compact sets whose (outer) boundary has the structure described in Remark 1. 
A preliminary version of the corresponding software along with the demos can be found at \cite{KSV24}.

Concerning interpolation, an appealing set is given by the so-called {\em Fekete points}, that are points which maximize the modulus of the Vandermonde determinant. As it is well-known, considering their Lebesgue constant  
it is immediate to get the (over)estimates  $\|L_n\|\leq n+1$ for the continuum Fekete points (since $\|\ell_j\|_{K}\leq 1$ for each $j$), and 
\begin{equation} \label{lebfek}
\|L_n\|\leq c_m(n+1)
\end{equation}
when they are obtained by extraction from a Chebyshev-like polynomial mesh $Z_n^m=\{z_1,\dots,z_M\}\subset \Gamma$. However, the continuum Fekete points are explicitly known only in two special instances, the interval (where they are the Gauss-Lobatto points) and the circle 
(where they are equally spaced in the arclength), in both cases with $\|L_n\|=\mathcal{O}(\log(n))$. 

On the other hand, the computation  
of Fekete points extracted from a polynomial mesh is known to be a NP-hard problem; cf. \cite{CMI09}.
Then, we can resort to points that approximately maximize the modulus of the Vandermonde determinant, extracting them from the polynomial mesh by  a greedy algorithm. Starting from \cite{BL08,SV09}, such {\em approximate Fekete points} have been computed by solving the underdetermined system 
\begin{equation} \label{afp}
\mathcal{V}_n^t \mathbf{u}=\mathbf{b}\;,\;\;
\mathcal{V}_n=V_n(Z_n^m)=[p_j(z_i)]\in \mathbb{C}^{M\times (n+1)} \;,\;\;Z_n^m=\{z_1,\dots,z_M\}\;,
\end{equation}
in a fixed polynomial basis $\{p_j\}$, where $\mathbf{b}$ is any nonzero vector, by $QR$ factorization with column pivoting. Indeed, the $n+1$ nonzero components of the solution vector $\mathbf{u}$ select the interpolation nodes $\Xi\subset Z_n^m$. This procedure corresponds to a greedy determinantal maximization, and the resulting interpolation points asymptotically behave as the continuum Fekete points, in the sense the corresponding uniform discrete probability measure converge weak-$\ast$ to the potential-theoretic equilibrium measure of the compact set; cf. \cite{BCLSV11,BL08} for a full discussion of these aspects, in both the univariate and the multivariate setting.

An interesting alternative is given by {\em Leja points}. For a fixed $\xi_1\in K$, the points are defined iteratively as $\xi_j=argmax_{z\in K}{\prod_{k=1}^{j-1}{|z-\xi_k|}}$, $j=2,\dots,n+1$, which means that differently from Fekete points they form a sequence, i.e. the first $\ell+1$ are Leja points for degree $\ell$. A relevant result has been recently proved by Totik \cite{T23}, who showed that the Lebesgue constant of Leja points has subexponential growth (a fact empirically well-known but missing before a theoretical base). On the other hand, it was previously proved in \cite{BBCL92} that Leja points behave asymptotically as Fekete points, in the sense the corresponding uniform discrete probability measure converge weak-$\ast$ to the potential-theoretic equilibrium measure of $K$. 

The same asymptotic property is shared by the two families of Leja-like sequences that we consider in this paper, namely {\em discrete Leja points} and {\em pseudo-Leja points}, both corresponding to a greedy discrete maximization on polynomial meshes. 
Indeed, discrete Leja points can be computed by $LU$ factorization with row pivoting of the whole matrix $\mathcal{V}_n$ in (\ref{afp}), where the pivots select the interpolation points within $Z_n^m$, a procedure substantially equivalent to compute $\xi_{j}=argmax_{z\in Z_n^m}{\prod_{k=1}^{j-1}{|z-\xi_k|}}$, $j=2,\dots,n+1$; cf. e.g. \cite{BDMSV10} for an analysis of this method, also in the multivariate setting. On the other hand, in our implementation we considered the pseudo-Leja points corresponding to the iteration $\xi_{j}=argmax_{z\in Z_{j-1}^m}{\prod_{k=1}^{j-1}{|z-\xi_k|}}$, $j=2,\dots,n+1$, after choosing the first point $\xi_{1}$ arbitrarily, e.g. $\xi_{1}$ is one of the points in $Z_1^m$ with larger imaginary component; cf. \cite{BCC12} (and \cite{K23} for a multivariate extension).

In Figures 2 and 4, we plot the Lebesgue constants 
of interpolation at approximate Fekete, discrete Leja and pseudo-Leja points, extracted from Chebyshev-like polynomial meshes $Z_n^m$ with $m=4$ on six different compact sets (see Figures 1 and 3), for degrees $n=1,\dots,50$. 
We also plot the Lebesgue constant of standard least-squares approximation of degree $n$ on the whole mesh. All these Lebesgue constants have been computed on the extraction meshes, recalling that with $m=4$ we get a relative error of less than $10\%$ and thus we substantially recover their actual size, cf. (\ref{asymptotic}). 

We make some observations on the main computational issues. In order to control the conditioning of the Vandermonde-like matrices, that becomes unacceptable already at moderate degrees with the standard monomial basis, we have chosen to work with a discrete orthonormalization of the shifted and scaled basis 
$q_j(z)=((z-z_b)/\delta)^{j-1}$, $1\leq j\leq n+1$, where $z_b=\frac{1}{M}\,\sum_{i=1}^ Mz_i$ is the barycenter of the mesh and $\delta=max_{1\leq i\leq M}|z_b-z_i|$. This approach allows to keep well-conditioning up to moderately high polynomial degrees (cf. e.g. \cite{BCLSV11}). Notice that, in view of (\ref{ls}) with $w_j\equiv1$, the Lebesgue constant on the mesh can then be simply computed via the relevant matrices as 
\begin{equation} \label{matrixleb}
\|\lambda_n\|_{Z_n^m}=\max_i\sum_j\left|\sum_k \pi_k(z_i)\overline{\pi_k(\xi_j)}\right|
=\|V_n(Z_n^m)R^{-1}Q^H\|_\infty\;,
\end{equation}
where $\{\xi_j\}$ are either the interpolation or the least-squares sampling points, and $V_n(\{\xi_j\})=QR$ the factorization of the corresponding Vandermonde-like matrix with $Q$ (rectangular) hermitian and $R$ square upper-triangular, that is $[\pi_1(z),\dots,\pi_{n+1}(z)]=[p_1(z),\dots,p_{n+1}(z)]R^{-1}$.

\begin{figure}[h!]
    \centering
        {\includegraphics[scale=0.2]{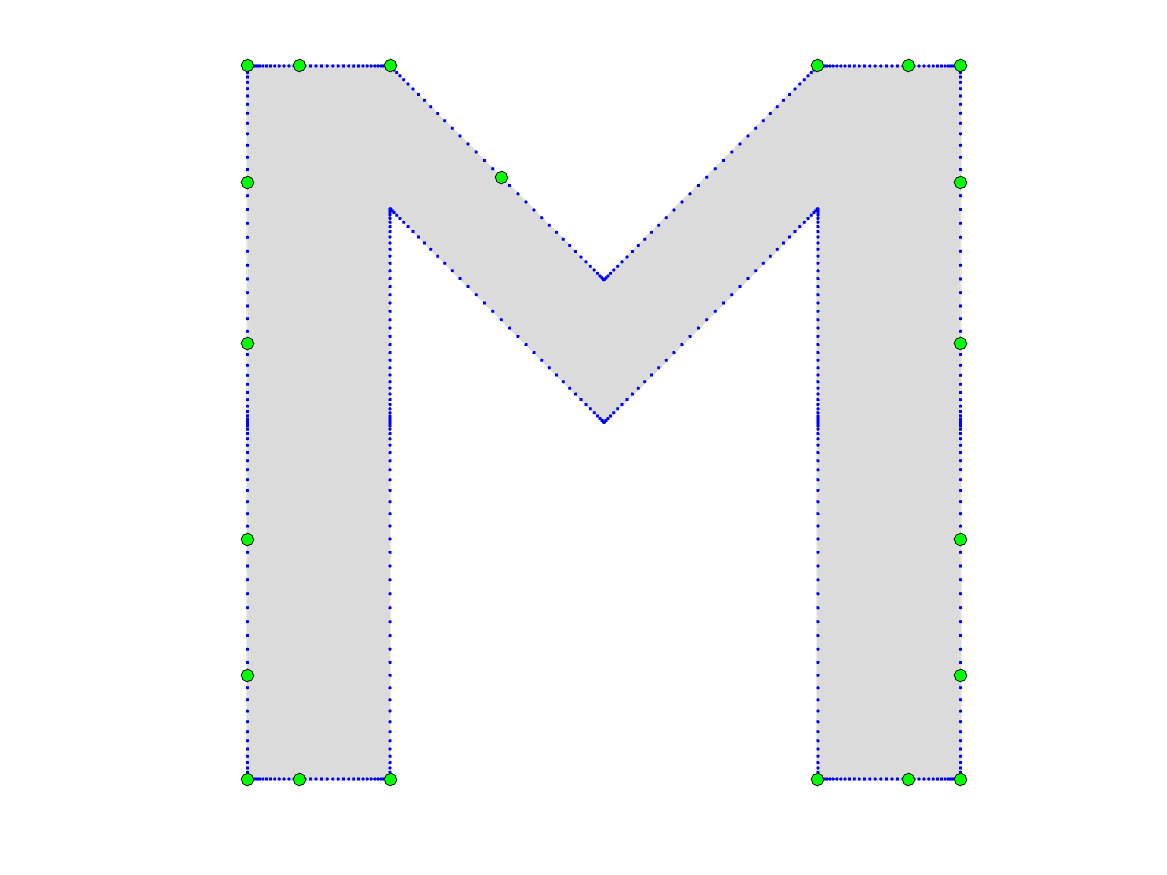}}
        {\hspace{-0.2cm}}
      {\includegraphics[scale=0.2]{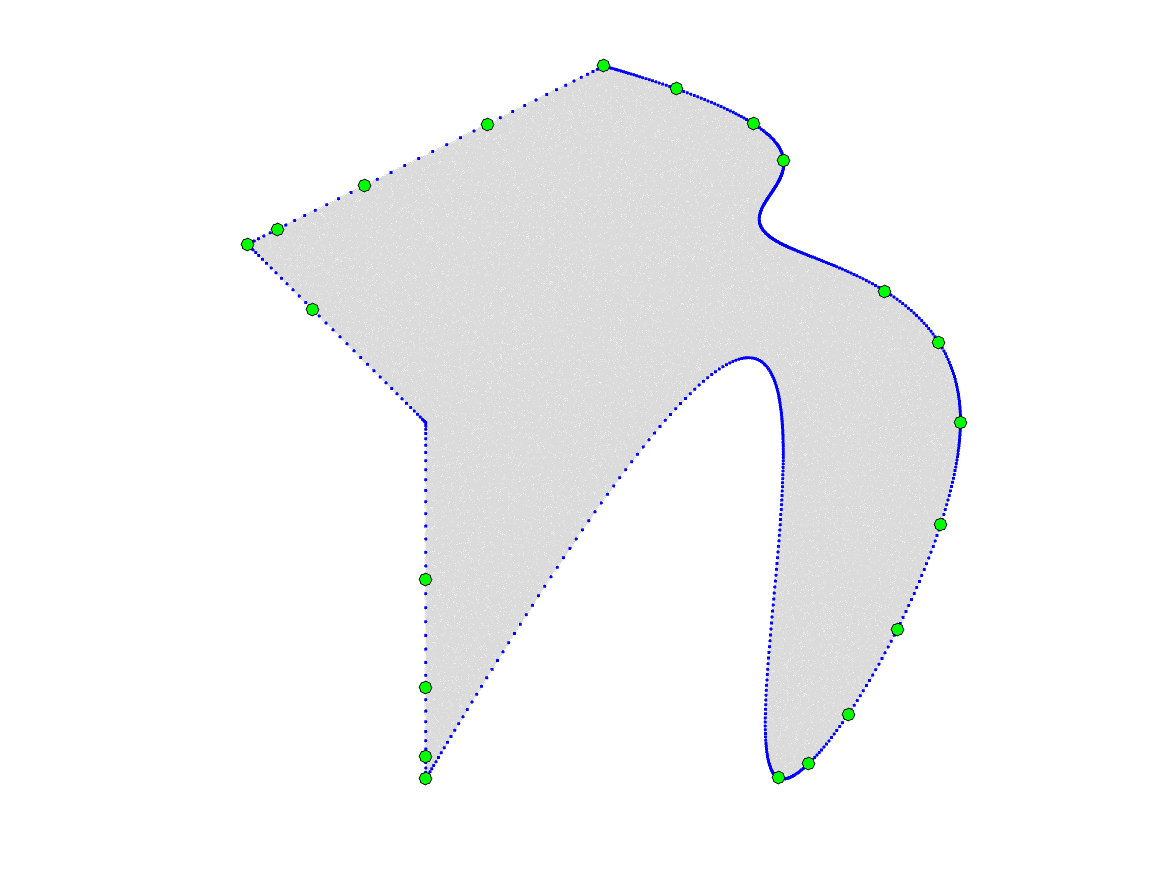}}
      {\hspace{-0.2cm}}
      {\includegraphics[scale=0.2] {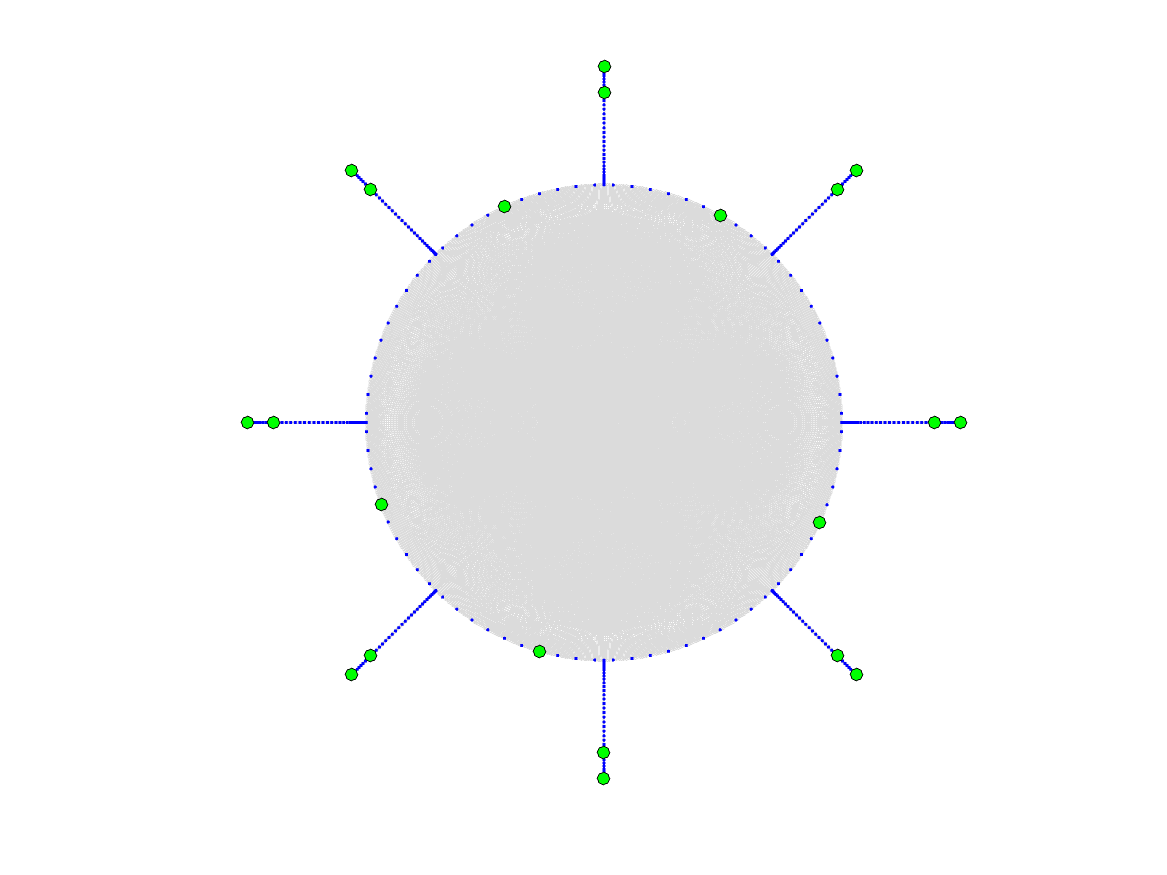}}
    \caption{\scriptsize{A polygon, a curvilinear polygon and a sun-shaped region as subsets of $\mathbb{C}$, the admissible polynomial mesh at degree 20 with  $m=2$ (blue dots), the 21 approximate Fekete points extracted from the mesh (green dots).}}
    \label{domains_1,2,3}
\end{figure}

\begin{figure}[h!]
    \centering
        {\includegraphics[scale=0.22]{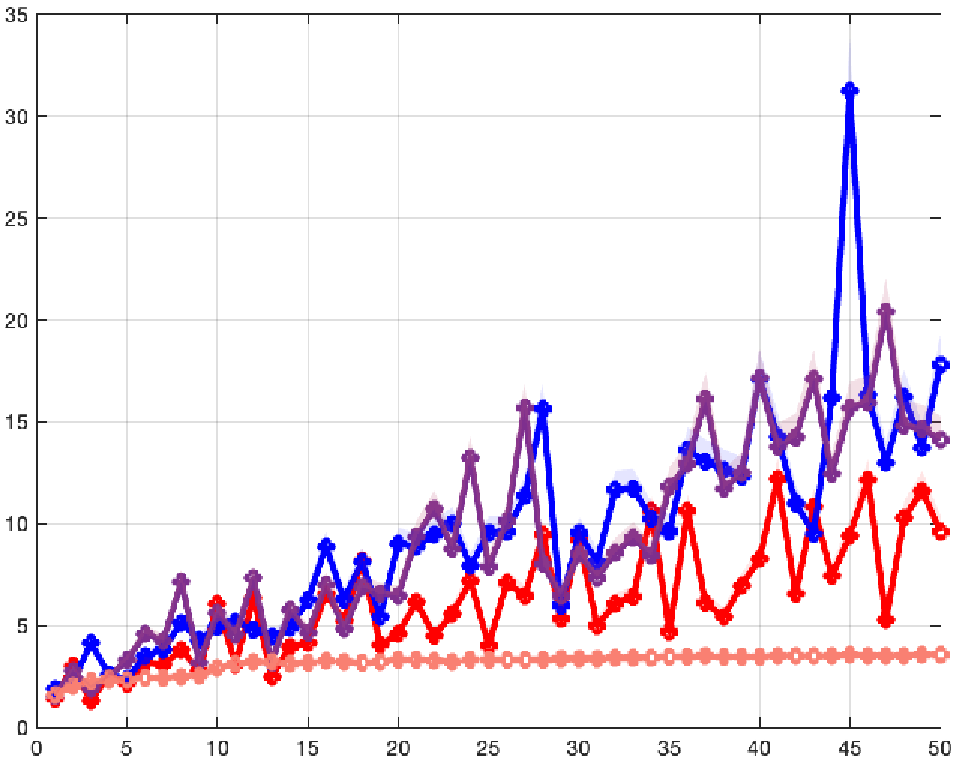}}
      {\hspace{0.2cm}}
        {\includegraphics[scale=0.22]{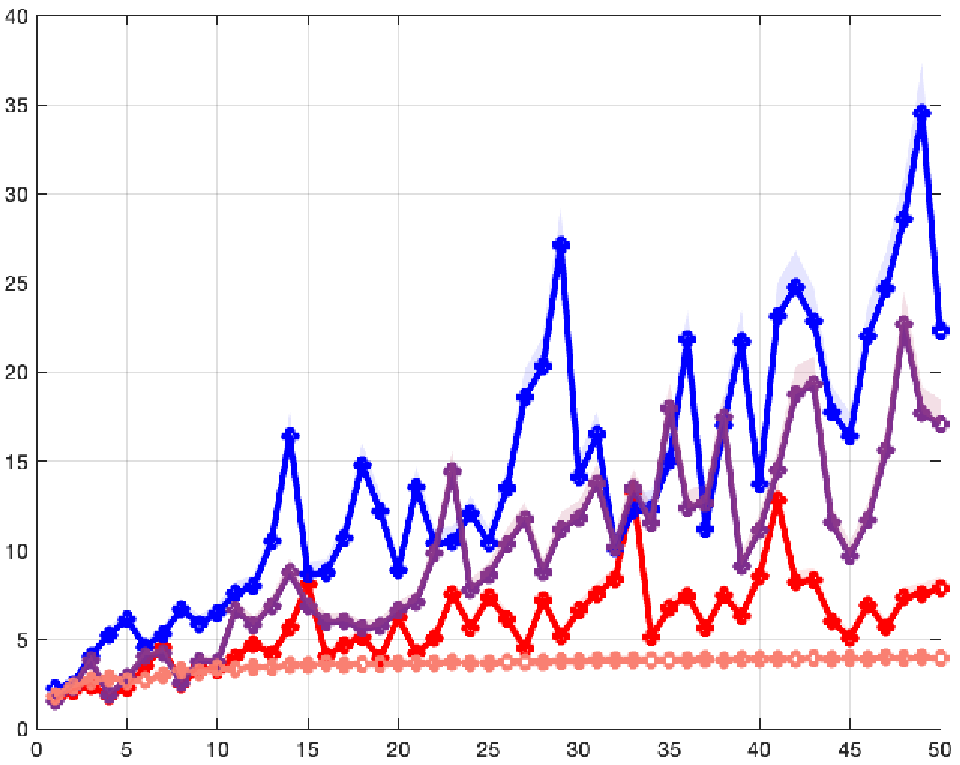}}
       {\hspace{0.2cm}}
        {\includegraphics[scale=0.22]{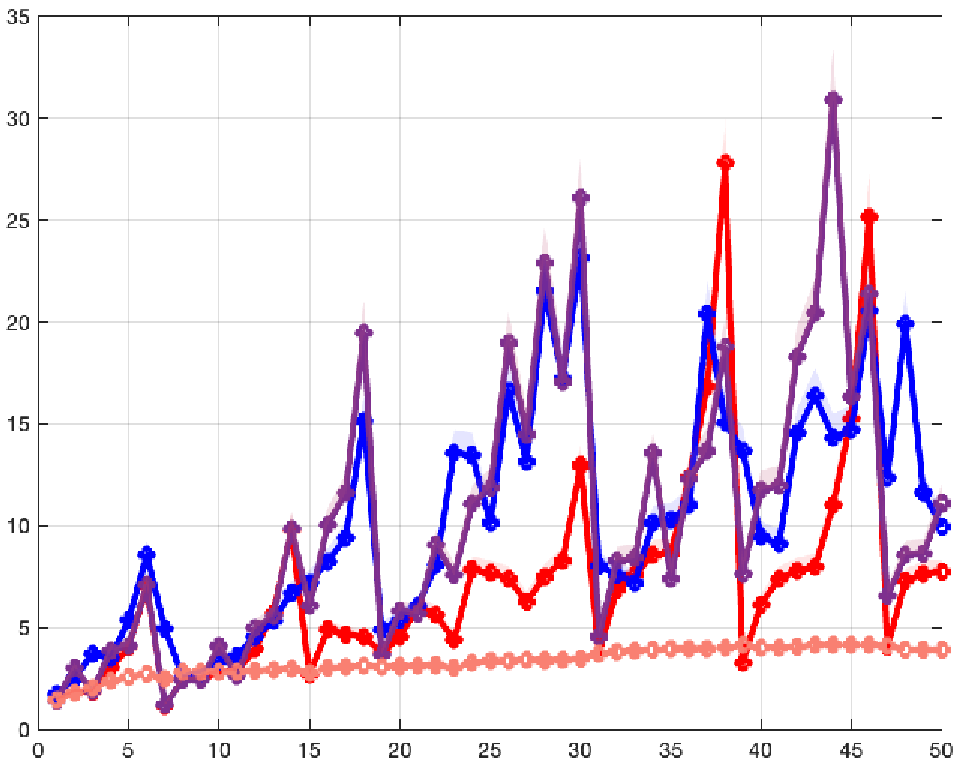}}
    \caption{\scriptsize{
    Lebesgue constants on the domains above for degrees $n=1,\dots,50$: least-squares on the whole mesh (pink dots), extracted approximate Fekete points (red dots), discrete Leja points (purple dots), pseudo-Leja points (blue dots). In these experiments, $m=4$.}}
    \label{leb_1,2,3}
\end{figure}

\begin{figure}[h!]
    \centering
        {\includegraphics[scale=0.2]{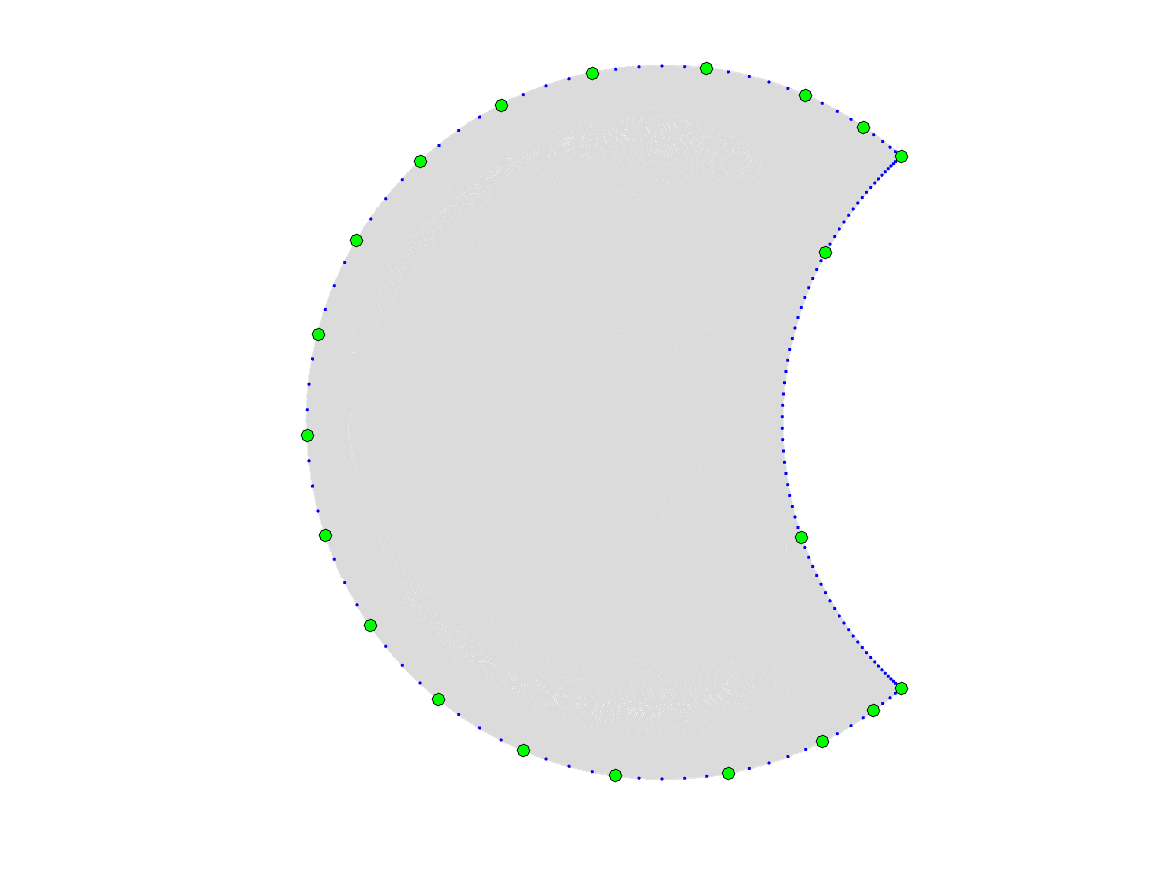}}
        {\hspace{-0.2cm}}
        {\includegraphics[scale=0.2]{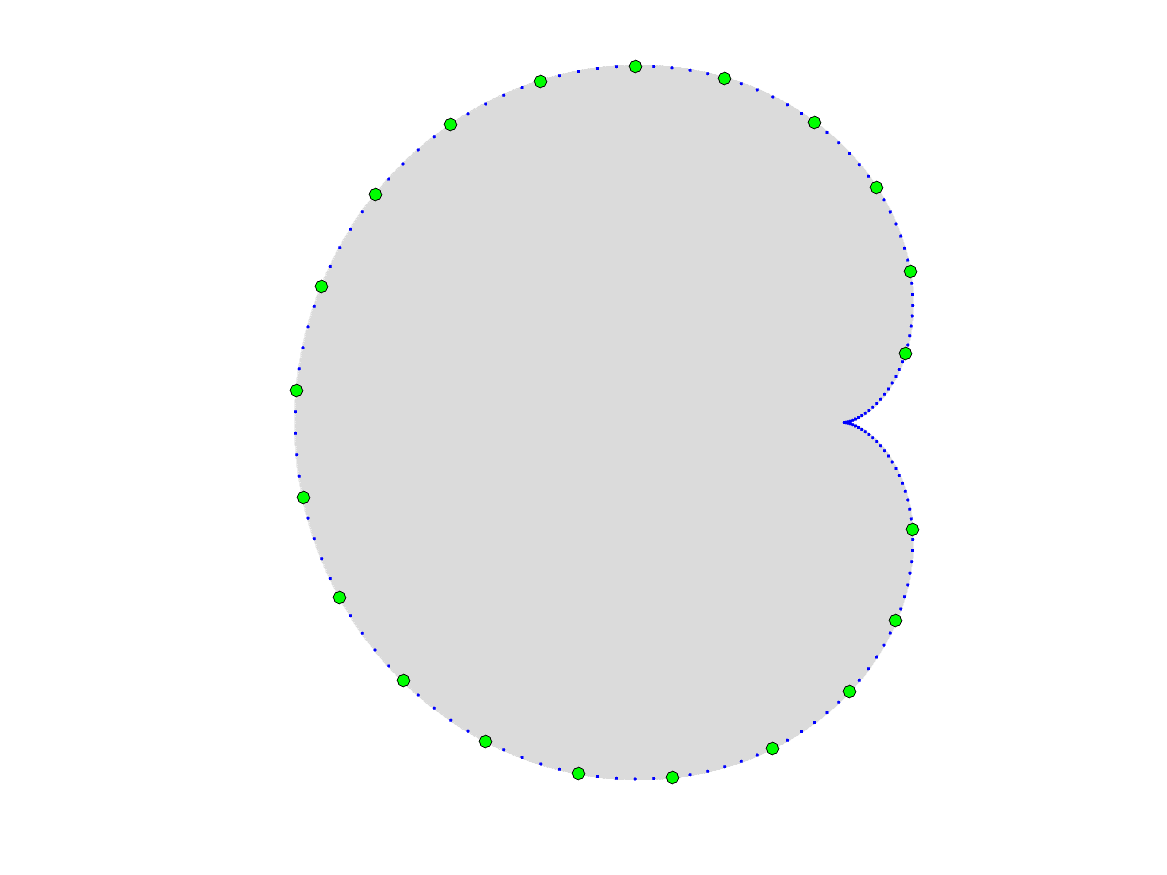}}
        {\hspace{-0.2cm}}
        {\includegraphics[scale=0.2]{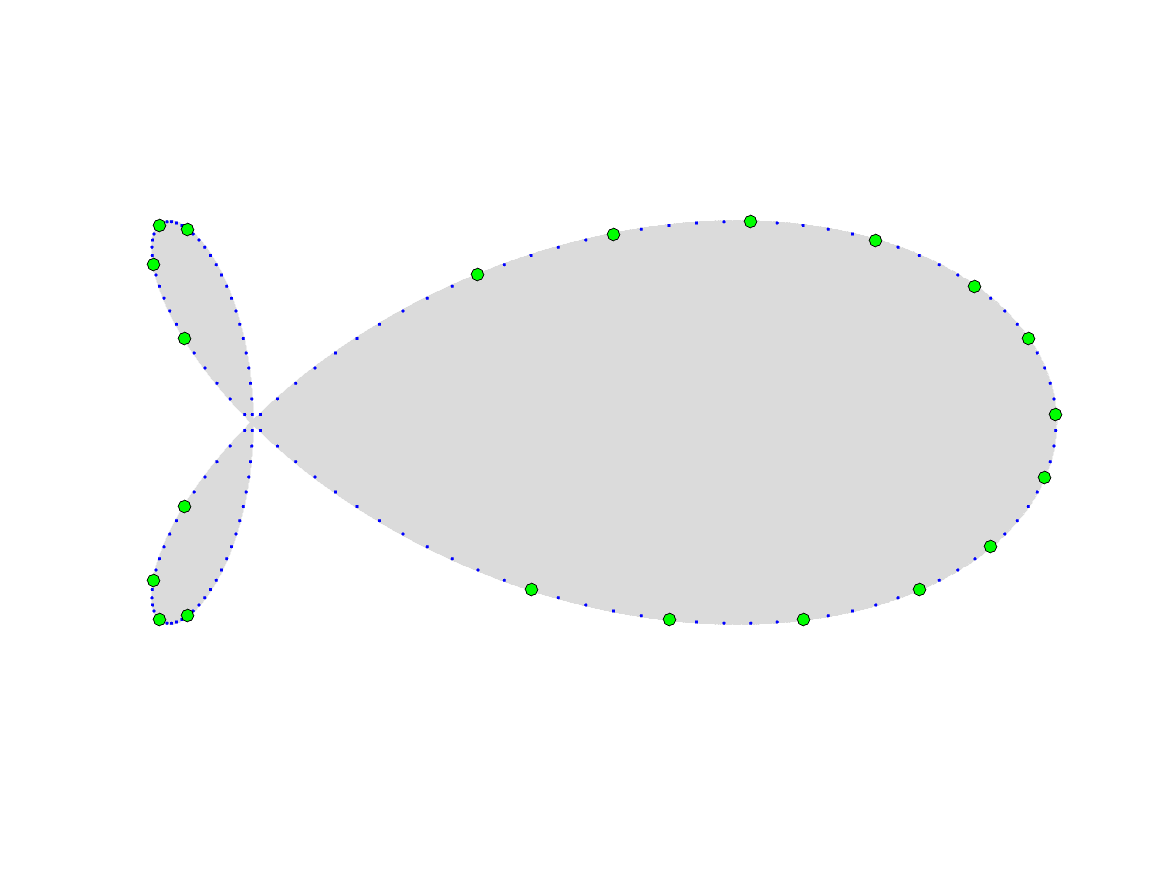}}
    \caption{\scriptsize{A lune, a cardioid and a torpedo as subsets of $\mathbb{C}$, the admissible polynomial mesh at degree 20 with  $m=2$ (blue dots), the 21 approximate Fekete points extracted from the mesh (green dots).}}
    \label{domains_4,5,6}
\end{figure}

\begin{figure}[h!]
    \centering
        {\includegraphics[scale=0.22]{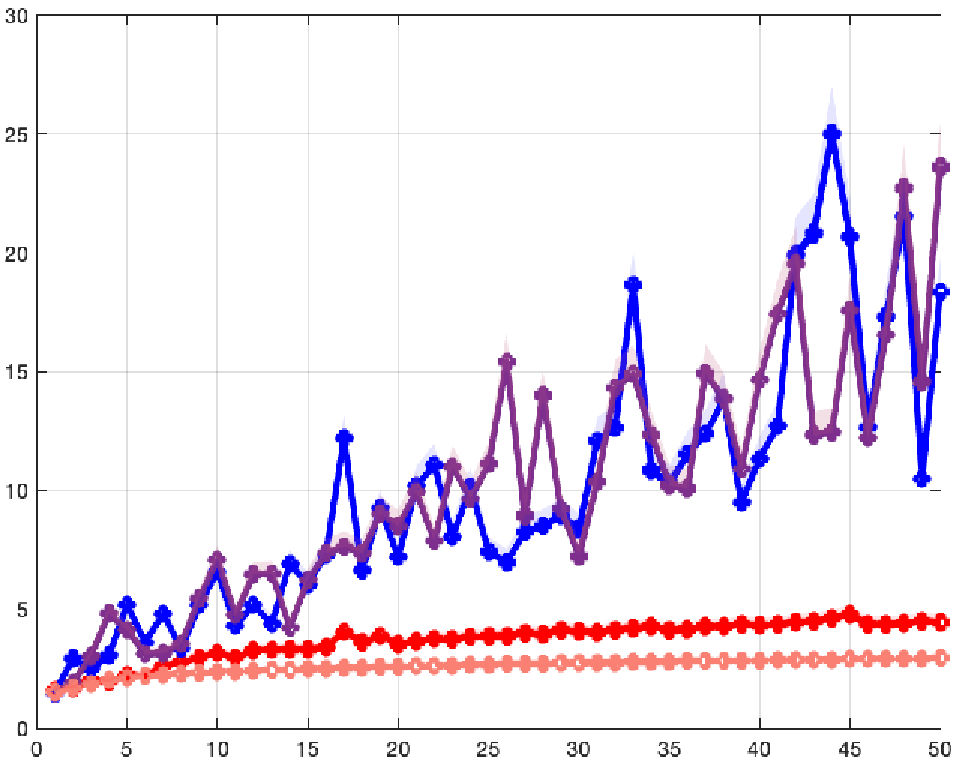}}
        {\hspace{0.2cm}}
        {\includegraphics[scale=0.22]{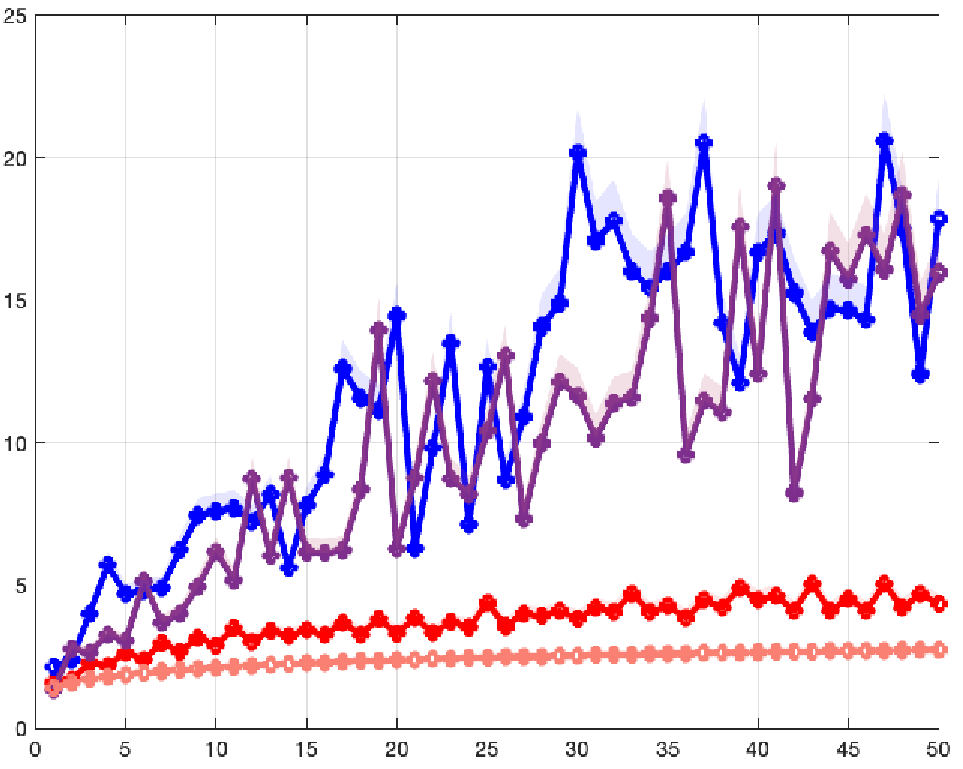}}
        {\hspace{0.2cm}}
        {\includegraphics[scale=0.22]{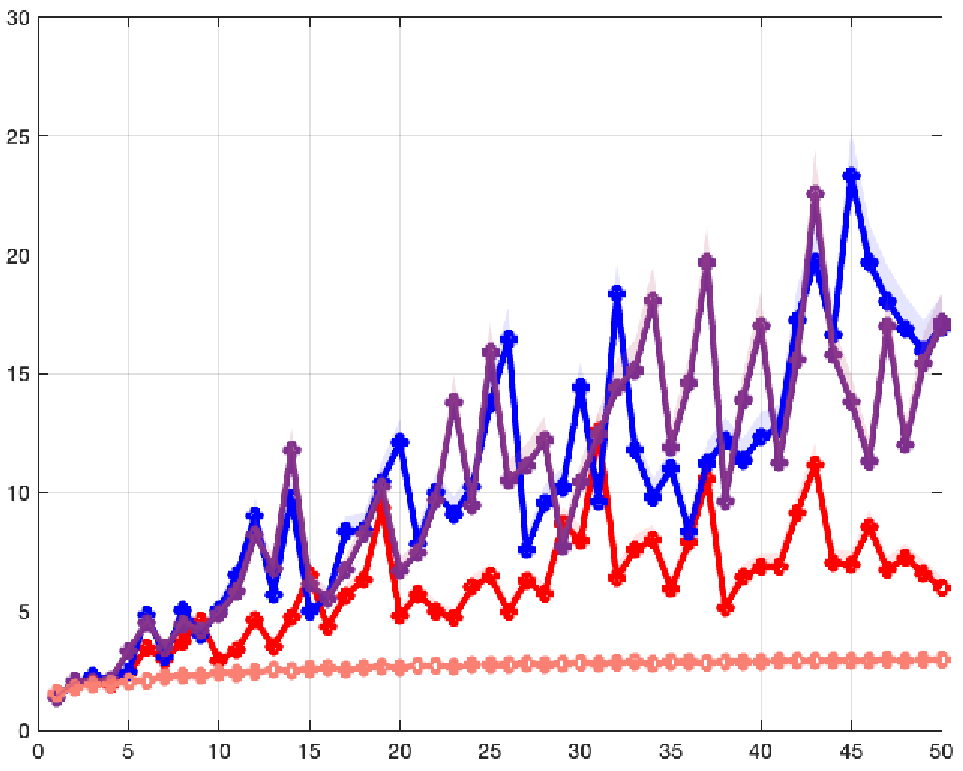}}
    \caption{\scriptsize{The Lebesgue constants as in Figure 2, on the three domains of Figure 3 (again, with $m=4$).}}
    \label{lebs_4,5,6}
\end{figure}

In the numerical experiments we have considered six complex regions whose (outer) boundaries can be tracked parametrically via one or more algebraic or trigonometric polynomials, in particular:
\begin{enumerate}
\item a M-shaped polygon with $12$ sides;
\item a curvilinear polygon, with  boundary  defined parametrically by linear and cubic splines;
\item a sun-shaped domain that consists of a unit disk and 8 rays that are segments of length $0.5$;
\item a lune defined as disk difference $B(-1,1.5)\setminus B(1,1.5)$;
\item  a cardioid, $\partial K$ being the closed curve
$$
z(t)=\cos(t)(1-\cos(t))+i(\sin(t)(1-\cos(t))), {\hspace{0.2cm}} t \in [0,2\pi];$$

\item a domain where $\partial K$ is the self-intersecting torpedo curve
$$
z(t)=\cos(t)\cos(2t)\exp(it), {\hspace{0.2cm}} t \in [0,2\pi].
$$
\end{enumerate}

First we observe that the interpolation points tend to privilege outward angles/tips/cusps as well as convex portions of the boundary and to avoid inward/concave portions, an electrostatic-like behavior that can be interpreted in connection with their potential theoretic background, cf. \cite{ST97}. We see that the Lebesgue constants of all discrete extremal sets show a slow increase, but those of Leja-like points have a more erratic behavior with larger  oscillations and tendentially higher values with respect to approximate Fekete points (a phenomenon already observed in the real multivariate setting, cf. e.g. \cite{BDMSV11}). On the other hand, Lebesgue constants of least-squares approximation on the whole polynomial mesh have the lowest values with an essentially logarithmic increase, staying below 5 up to degree 50 in all the six examples.

\vskip1cm 
\noindent
{\bf Acknowledgements.} 
Work partially supported by the DOR funds of the University of Padova and by the INdAM-GNCS (A. Sommariva, M. Vianello), and 
by the National Science Center - Poland, grant Preludium Bis 1, N. 2019/35/O\-/ST1/02245 (D.J. Kenne). The research cooperation was funded by the program Excellence Initiative – Research University at the Jagiellonian University in Krak\`{o}w (A. Sommariva). 
This research has been accomplished within the RITA ``Research ITalian network on Approximation" and the SIMAI Activity Group ANA\&A (A. Sommariva, M. Vianello), and the UMI Group TAA ``Approximation Theory and Applications" (A. Sommariva).

\end{document}